\documentclass[11pt]{article}

\usepackage{graphicx}
\usepackage{amssymb}
\newtheorem{definition}{Definition}

\newtheorem{theorem}{Theorem}
\newtheorem{lemma}{Lemma}

\newcommand{\EQ}{\begin{equation}\begin{array}{lllllllll}}
\newcommand{\EE}{\end{array}\end{equation}}
\newcommand{\EQQ}{$$\begin{array}{lllllllll}}
\newcommand{\EEE}{\end{array}$$}
\newcommand{\MT}{\left[ \begin{array}{ccccccccc}}
\newcommand{\EM}{\end{array}\right]}
\newcommand{\bean}{\begin{equation}\begin{array}{rcllllllll}}
\newcommand{\eean}{\end{array}\end{equation}}
\newcommand{\bea}{$$\begin{array}{rcllllllll}}
\newcommand{\eea}{\end{array}$$}
\def\ds{\displaystyle}
\def\Fr{\ds \frac}
\def\=={&=&}
\def\dsum{\ds\sum}
\def\Real{{I\!\!R}}

\def\calH{{\cal H}}

\def\L2O{X}
\title{\LARGE \bf
Partial Observability and its Consistency for PDEs
\thanks{This work was supported in part by NRL and AFOSR. The views expressed in this document are those of the author and do not reflect the official policy or position of the Department of Defense or the U.S. Government.}
}

\author{Wei Kang\thanks{Wei Kang is with Faculty of Applied Mathematics, Naval Postgraduate School, Monterey, CA, USA, {\tt\small wkang@nps.edu}}, Liang Xu\thanks{Liang Xu is at Naval Research Laboratory, Monterey, CA, USA, {\tt\small liang.xu@nrlmry.navy.mil}}, and Francis X. Giraldo\thanks{Francis X. Giraldo is with Faculty of Applied Mathematics, Naval Postgraduate School, Monterey, CA, USA,
        {\tt\small fxgirald@nps.edu}}
}

\begin{document}

\maketitle
\thispagestyle{empty}

\begin{abstract}
In this paper, a quantitative measure of partial observability is defined for PDEs. The quantity is proved to be consistent if the PDE is approximated using well-posed approximation schemes. A first order approximation of an unobservability index using an empirical Gramian is introduced. Several examples are presented to illustrate the concept of partial observability, including Burgers' equation and a one-dimensional nonlinear shallow water equation. 

\end{abstract}

\section{Introduction}
Observability is a fundamental property of dynamical systems \cite{isidori,kailath} with an extensive literature. It can be considered as a measure of well-posedness for the estimation of system states using sensor information as well as additional user knowledge about the system. We do not give a review of the huge literature here on this subject. Some interesting work can be found in \cite{gauthier,hermann,xia} for nonlinear systems, \cite{zuazua} for partial differential equations (PDEs), \cite{mohler} for stochastic systems, and \cite{zheng} for normal forms. 

For some models of high dimensional systems, the traditional concept of observability is not applicable. For instance, models used in numerical weather prediction have millions of state variables. Some variables are strongly observable while some others are extremely weakly observable. It is known that the sparse sensor network cannot make the entire system uniformly observable. These types of problems call for a {\it partial observability} analysis in which we study the observability of finite many modes (or state variables) that are important for weather prediction and ignore the modes that are less important. In addition, it is desirable to measure the observability of a system {\it quantitatively}. It is not good enough to just tell that a set of variables is observable or not. It is important to tell how strong or weak the observability is. In the large amount of data collected for weather prediction, only $5\%$ or less are actually useful for each prediction. Finding a high-value data set that improves the level of observability requires a quantitative measure of observability. When moving sensors are used, a quantitative measure of observability is fundamental in finding  optimal sensor locations. Another issue about models of high dimensional systems is how to practically verify their observability when a model is given as a numerical input-output function, such as a code, rather than a set of differential equations. The definition of observability should be {\it computational} so that one can numerically verify the concept.  

In \cite{kang-xu1,kang-xu2}, a definition of observability is introduced using dynamic optimization. This concept  is developed in a project of finding the best sensor locations for data assimilations, a computational algorithm widely used in numerical weather prediction. Different from traditional definitions of observability, the one in \cite{kang-xu1,kang-xu2} is able to collectively address several issues in a unified framework, including a quantitative measure of observability, partial observability, and improving observability by using user knowledge. Moreover, computational methods of dynamic optimization provide practical tools of numerically approximating the observability of complicated systems that cannot be treated using an analytic approach. 

In this paper we extend the definition of observability in \cite{kang-xu1,kang-xu2} to systems defined by PDEs. A quantitative measure of partial observability makes perfect sense for infinite dimensional systems such as PDEs. However, its computation is carried out using finite dimensional approximations. It is known in the literature that an approximation of a PDE using ordinary differential equations (ODEs) may not preserve the property of observability, even if the approximation scheme is convergent and stable \cite{zuazua,cohn}. Therefore, to develop the concept of partial observability for PDEs, it is important to understand its consistency in ODE approximations. In Section \ref{sec2}, some examples from the literature are introduced to illustrate the issues being addressed in this paper. Observability is defined for PDEs in Section \ref{sec3}. In Section \ref{sec4}, a theorem on the consistency of observability is proved. The relationship between the unobservability index and an empirical Gramian is addressed in Section \ref{sec5}, which serves as a first order approximation of the observability. The theory is verified using several examples in Section \ref{sec7}.

\section{Some issues on observability}
\label{sec2}
The theory in this paper is developed for problems that are so large and complicated that the traditional method of analysis does not apply. Before we introduce the concept and theorems,  we first use a few simple examples to illustration some issues to be addressed. Consider the initial value problem of a heat equation
$$\begin{array}{lll}
u_t(x,t)=u_{xx}(x,t), \; x\in [0, L], t\in [0, T]\\
u(0,t)=u(L, t)=0\\
u(x,0)=f(x)
\end{array}$$
Suppose the measured output is 
$$y(t)=u(x_0,t)$$
for some $x_0\in [0,L]$. In this example, we assume $L=2\pi$, $T=10$, and $x_0=0.5$. The solution and its output have the following form
$$\begin{array}{lll}
u(x,t)=\ds\sum_{k=1}^{\infty} \bar u_k(t)\sin \left(\Fr{k\pi x}{L}\right)\\
y(t)=\ds\sum_{k=1}^{\infty} \bar u_k(t)\sin \left(\Fr{k\pi x_0}{L}\right)
\end{array}$$
where the Fourier coefficients satisfy an ODE
$$\Fr{ d\bar u_k}{dt}=-\left(\Fr{k\pi}{L}\right)^2 \bar u_k$$
Define 
$$\begin{array}{lll}
u^N=\MT \bar u_1 & \bar u_2 & \cdots \bar u_N\EM^T\\
A^N=\mbox{diag}\left(\MT \left(\Fr{\pi}{L}\right)^2 & \left(\Fr{2\pi}{L}\right)^2&\cdots &\left(\Fr{N\pi}{L}\right)^2\EM\right)\\
C^N=\MT \sin\left(\Fr{\pi x_0}{L}\right) & \sin\left(\Fr{2\pi x_0}{L}\right)&\cdots &\sin\left(\Fr{N\pi x_0}{L}\right)\EM
\end{array} $$
A $N$th order approximation of the original initial value problem with output is defined by a system of ODEs
\EQ
\label{eqheatode}
du^N/dt=-A^Nu^N\\
u^N(0)=u^N_0\\
y=C^Nu^N
\EE
A Gramian matrix \cite{kailath} can be used to measure the observability of $u^N_0$. More specifically, given $N>0$ the observability Gramian is
$$W=\ds\int_0^Te^{(-A^N)'t}(C^N)'C^Ne^{-A^Nt}dt$$
Its smallest eigenvalue, $\sigma^N_{min}$, measures the observability of $u^N_0$. A small value of $\sigma^N_{min}$ implies weak observability. If the maximum sensor error is $\epsilon$, then the worst estimation error of $u^N_0$ is bounded by 
\EQ
\label{esterror}
\Fr{\epsilon}{\sqrt{\sigma^N_{min}}}
\EE

\begin{figure}[!ht]
	\begin{center}
		\includegraphics[width=3.0in,height=2.0in]{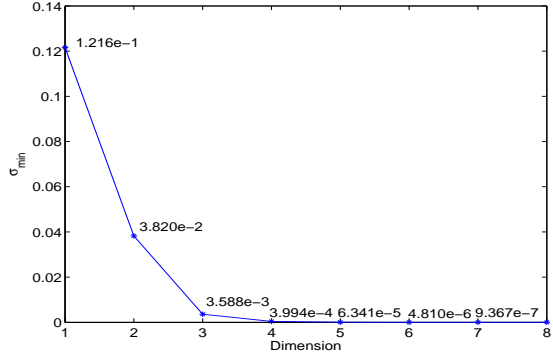} 
		\caption{Observability of heat equation}
		\label{Figheat}
		\end{center}
\end{figure}

The system has infinitely many modes in its Fourier expansion. However, it has a single output. The computation shows that the output can make the first mode observable. However, when the number of modes is increased, their observability decreases rapidly. From Figure \ref{Figheat}, for $N=1$ we have $\sigma^N_{min}=1.216\times 10^{-1}$, which implies a reasonably observable $\bar u_1(0)$. However, when $N$ is increased, the observability decreases rapidly. For $N=8$, $\sigma^N_{min}$ is almost zero, i.e 
$$\MT \bar u_1(0)&\bar u_2(0)&\cdots &\bar u_8(0)\EM^T$$ 
is extremely weakly observable, or practically unobservable. According to (\ref{esterror}), a small sensor noise results in a big estimation error.  For this problem, it is not important to achieve observability for the entire system. All we need is the partial observability for the critical modes. 

For a simple system like the heat equation, we already know that the high modes do not affect the solution. Their observability is not important. One may choose a discretization scheme based on the important modes only so that the observability is achieved using a simplified model. However, for large and highly nonlinear systems with tens of thousands or even millions of state variables, such as in a process of numerical weather prediction or power network control, model reduction or changing the scheme of discretization is almost impossible because the models are already packaged in the form of software and copyright issues are likely involved. For these types of applications, a {\it quantitative measure} of {\it partial observability} is useful for several reasons. If a finite number of modes is enough to achieve accurate state approximates, guarantee the observability of these finite modes is a practical solution. In large-scale networded systems with a very high dimension, operators may focus on a local area at a given period of time. In this case, the observability of the entire system is irrelevant. It is useful to achieve a partial observability just for the states directly related to the area of focus. For moving sensors, a quantitative measure of observability can be used as a cost function in finding optimal sensor locations in which the observability is maximized. 

Another issue to be addressed in this paper is {\it consistency}. In general the observability for PDEs is numerically computed using a system of ODEs as an approximation. However, it is not automatically guaranteed that the observability of the ODEs is consistent with the observability of the original PDE. In fact, a convergent discretization of a PDE may not preserve its observability. Take the following wave equation as an example
\EQ
\label{wave}
u_{tt}-u_{xx}=0, &0<x<L, \; 0<t<T\\
u(0,t)=u(L,t)=0, & 0<t<T\\
u(x,0)=u_0(x), u_t(x,0)=u_1(x), &0<x<L
\EE
It is known that the total energy of the system can be estimated by using the energy concentrated on the boundary. However, in \cite{zuazua} it is proved that the discretized ODEs do not have the same observability. The energy of solutions is given by
$$E(t)=\Fr{1}{2}\ds\int_0^L \left( |u_t(x,t)|^2+|u_x(x,t)|^2  \right) dx$$
This quantity is conserved along time. It is known that, when $T>2L$, the total energy of solutions
can be estimated uniformly by means of the energy concentrated on the boundary $x=L$. More precisely, there exists $C(T)>0$ such that 
\EQ
\label{waveobs}
E(0) \leq C(T)\ds\int_0^T |u_x(L,t)|^2dt
\EE
Now consider the discretized system using a finite difference method,
\EQ
\label{wave2}
u''_j(t)=\Fr{u_{j+1}(t)+u_{j-1}(t)-2u_j(t)}{h^2}, & 0<t<T, \;j=1,2,\cdots,N\\
u_0(t)=u_{N+1}(t)=0, & 0<t<T\\
u_j(0)=u^0_j, u'_j(0)=u^1_j, &j=0,1,\cdots,N+1
\EE
The total energy of the ODEs is given by
$$E_h(t)=\Fr{h}{2}\ds\sum_{j=0}^N \left(|u'_j(t)|^2+\left|\Fr{u_{j+1}(t)-u_j(t)}{h}   \right|^2   \right)$$
This quantity is conserved along the trajectories of the ODEs. The energy on the boundary is defined by
$$\ds\int_0^T \left| \Fr{u_N(t)}{h} \right|^2dt$$
Because the solutions of (\ref{wave2}) converges to the solutions of (\ref{wave}), we expect that the total energy of (\ref{wave2}) can be uniformly estimated using the energy concentrated along a boundary, i.e. the following inequality similar to (\ref{waveobs}) holds for some $C(T)$ uniformly in $h$,
$$E_h(0)\leq C(T)\ds\int_0^T \left| \Fr{u_N(t)}{h} \right|^2dt$$
However, it is proved in \cite{zuazua} that this inequality is not true. In fact, the ratio between the total energy and the energy along the boundary is unbounded as $h\rightarrow 0$. To summarize, the observability of a PDE is not necessarily preserved in its discretizations. 

In this paper, we introduce a {\it quantitative measure} of {\it partial observability} for PDEs. Sufficient conditions are proved for the {\it consistency} of the observability for well-posed discretization schemes. 

\section{Partial Observability}
\label{sec3}
Consider a nonlinear initial value problem
\EQ
\label{eqpdemodelnon}
u_t=F(t,u,u_x,\cdots), &\mbox{in } \Omega \times (0, T]\\
u=u_0 & \mbox{in } \Omega \times \{ t=0\}\\
y_u(t)=\calH (u(\cdot,t))
\EE
where $\Omega$ is an open set in $\Real^n$, $F$ is a continuous function of $t$, $u$, and its derivatives with respect to $x$. Let $X$ be a Banach space of functions defined on $\Omega$. The initial condition, $u_0$, lies in a subspace, $ V_0$, of $X$. In the following, $u(t)$ represents $u(\cdot,t)$. A solution $u(t)$ of (\ref{eqpdemodelnon}), in either strict or weak sense, is a $\L2O$-valued function in a subspace $V$ of $C^0([0,T],X)$. If additional boundary conditions are required, we assume that all functions in $V$ satisfy the  boundary conditions. We assume that (\ref{eqpdemodelnon}) is locally {\it well-posed} in the Hadamard sense (\cite{hille,richtmyer}) around a nominal trajectory. More specifically, let $u(t)$ be an nominal trajectory. We assume that there is an open neighborhood $D_0\subset V_0$ that contains $u(0)$, such that  
\begin{itemize}
\item For any $u_0\in D_0$, (\ref{eqpdemodelnon}) has a solution.
\item The solution is unique.
\item The solution depends continuously on its initial value.
\end{itemize}
Proving the well-posedness of nonlinear PDEs is not easy. Nevertheless, for well-posed problems the consistency of observability is guaranteed. 

In (\ref{eqpdemodelnon}), $\calH (u(\cdot,t))$ or in short notation $y_u(t)$, represents the output of the system associated with the solution $u(x,t)$, where $\calH$ is a mapping, linear or nonlinear, from $\L2O$ to $\Real^p$. We assume that $y_u(\cdot)$ stays in a normed space of functions from $[t_0,t_1]$ to $\Real^p$. Its norm is denoted by $|| \cdot ||_Y$. We say that $\calH $ is continuous in a subset of $C^0([0,T],X)$ if
for any sequence
$\{u_k(t)\}_{k=k_0}^\infty$ in the subset and  function $u(t)\in V$,  
\EQ
\label{continuityassumption}
u_k(t) \rightarrow u(t) \mbox{ uniformly on } [0,T] \mbox{ implies } \ds\lim_{k\rightarrow \infty} ||y_{u_k}-y_u||_Y=0 
\EE

Instead of the entire state space, the observability is defined in a finite dimensional subspace. Let 
$$W=\mbox{span}\{ e_1,e_2,\cdots,e_s\}$$
be a subspace of $V_0$ generated by a basis $\{ e_1,e_2,\cdots,e_s\}$. In the following, we analyze the observability of a component of $u(0)$ using estimates from $W$. Therefore $W$ is called the {\it space for estimation}.

Let $u(t)$ be a nominal trajectory satisfying (\ref{eqpdemodelnon}). Suppose $W\bigcap D_0\neq \emptyset$ and $u(0)$ has a best estimate in $W\bigcap D_0$, denoted by $u_w(0)$, in the sense that $u_w(t)$ minimizes the following output error,
\EQ
\label{uproj}
\min||y_{u_w}-y_u||_Y\\
\mbox{subject to} \\
du_w/dt=F(t,u_w,\cdots), &\mbox{in } \Omega \times (0, T]\\
u_w(0)\in W\bigcap D_0
\EE
Let $u_r(t)=u(t)-u_w(t)$ be the remainder, then 
\EQ
\label{decomposition}
u(t)=u_w(t)+u_r(t)
\EE
If the output $y_u(t)$ represents the sensor measurement, then it has noise. The data that we use in a estimation process has the following form 
$$y_u(t)+d(t)$$
where $d(t)$ is the measurement error. The observability addressed in this paper is a quantity that defines the sensitivity of the estimation error relative to $d(t)$. From (\ref{decomposition}), the best estimate $u_w(t)$ has an error that is the remainder $||u_r(0)||_X$. This error is not caused by $d(t)$ because the remainder cannot be reduced no matter how accurate the output is measured. This error is due to the choice of $W$, not the observability of $W$. Therefore, the observability is defined for $u_w(0)$ only, thus a {\it partial observability}. For a strongly observable $u_w(0)$, its estimate may not be close to $u(0)$ if the remainder, $u_r(0)$, is large. If the goal is to estimate $u(0)$, an observable $u_w(0)$ is useful only if the value of $u_r(0)$ is either known or small. In the rest of the paper, we assume $u_r(0)=0$ and $u(0)\in W$. In applications, the concept is applicable to $u(0)$ with a small $u_r(0)$, (see, for instance, \cite{kang-xu3}).

\begin{definition}
\label{def1}
Given a nominal trajectory $u$ of (\ref{eqpdemodelnon}). Suppose 
$$u(0) \in W\bigcap D_0.$$ For a given $\rho >0$, suppose the sphere in $W$ of radius $\rho$ centered at $u(0)$ is contained in $D_0$. We define
\EQ
\label{eqcostnon}
\epsilon = \inf ||y_{\hat u}-y_u||_Y
\EE
where $\hat u$ satisfies
\EQ
\label{hatunon}
\hat u_t=F(t,\hat u,\hat u_x,\cdots)\\
\hat u(0)\in W\bigcap D_0\\
||\hat u(0) -u(0)||_X=\rho 
\EE
Then $\rho/\epsilon$ is called the unobservability index of $u(0)$ along the trajectory $u(t)$ at the level of $\rho$.
\end{definition}

\noindent{\bf Remark}. The ratio $\rho/\epsilon$ can be interpreted as follows: if the maximum error of the measured output, or sensor error, is $\epsilon$, then the worst estimation error of $u(0)$ is $\rho$. Therefore, a small value of $\rho/\epsilon$ implies strong observability of $u(0)$. Different from most traditional definitions of observability, $W$ is a subspace of the state space. However, in the special case that $W$ equals the entire state space that has a finite dimension and if the system is linear,  $\epsilon^2/\rho^2$ equals the smallest eigenvalue of the observability Gramiam (See Section \ref{sec5}). \\
\\
{\bf Remark}. The effectiveness of this concept is verified in \cite{kang-xu3} using Burgers' equation. A 4D-Var data assimilation, a method of state estimation in weather prediction, is applied to several data sets from different sensor locations. The result shows that the estimation using sensor locations with a smaller unobservability index results in more accurate estimates than those using the data from other  sensor locations with higher unobservability indices. 
\\

To numerically compute a system's observability, (\ref{eqpdemodelnon}) is approximated using ODEs. In this paper, we consider a general approximation scheme using a sequence of ODEs,
\EQ
\label{eqodemodelnon}
du^N/dt=F^N(t,u^N), & u^N\in \Real^N\\
u^N(0)=u^N_0
\EE
and two sequences of linear mappings
\EQ
\label{mappingnon}
P^N: V_0 \rightarrow \Real^N\\
\Phi^N: \Real^N \rightarrow \L2O
\EE
The norm in $\Real^N$ is represented by $||\cdot ||_N$. We assume that the approximation scheme is  well-posed, more specifically
\begin{itemize}
\item For any bounded set $B$ in $D_0$, there exist $M>0$ and $\alpha >0$ so that any solution of (\ref{eqpdemodelnon}) with $u(0)\in B$ satisfies
\EQ
\label{rateconv}
||u(t)-\Phi^N(u^N(t))||_X \leq \Fr{M}{N^\alpha}
\EE
uniformly in $[0, T]$, where $P^N(u(0))=u^N(0)$. 
\item For any $u\in W$,
\EQ
\label{mapping3b}
||u||_X=||P^N(u)||_N+a_N||P^N(u)||_N\\
\ds\lim_{N\rightarrow \infty}a_N=0
\EE
\end{itemize}

Given the space for estimation $W$, we define a sequence of subspaces, $W^N \subseteq\Real^N$, by
$$W^N=P^N(W)$$
They are used as the space for estimation in $\Real^N$. If $\{e_1,e_2,\cdots,e_s\}$ is a basis of $W$, then their projections to $W^N$ are denoted by
$$e_i^N=P^N(e_i), \;\; i=1,2,\cdots,s$$  
So $W^N=\mbox{span}\{ e_1^N, e_2^N,\cdots,e_s^N\}$.\\

\noindent{\bf Example}.
For a spectral method, approximate solutions can be expressed in terms of an orthonormal basis 
$$\{ q_k(x): k=0,1,2,\cdots \}$$
For any function, 
$$v(x)=\sum_{k=0}^\infty v_kq_k(x)$$ one can define 
\EQ
\label{eq_spectral1}
P^N (v)=\MT v_0&v_1&\cdots&v_N\EM^T
\EE
Obviously, $\Phi^N$ is defined by
\EQ
\label{eq_spectral2}
\Phi^N(\MT v_0&v_1&\cdots&v_N\EM^T)=\dsum_{k=0}^N v_kq_k
\EE
In the case of $l^2$-norm for all $\Real^N$, $||\cdot||_N$ is consistent with $||\cdot||_X$ if $X=L^2(\Omega)$. Typically, the space for estimation consists of finite many modes
\EQ
\label{ad_example_W}
W=\mbox{span} \{ q_1(x), q_2(x), \cdots, q_s(x)\}
\EE
$\square$\\

\noindent{\bf Example}. Some approximation methods, such as finite difference and finite element, are based on a grid defined by a set of points in space, 
$$\{ x_k^N\}_{k=0}^N$$
and a basis $\{ q_k^N \}$ satisfying
\EQ
\label{eq_grid1}
q_k^N(x_j^N)=\left\{ \begin{array}{lll} 1 & k=j\\
0, & otherwise
\end{array}\right.
\EE
In this case, the mappings in the approximation scheme is defined as follows
\EQ
\label{eq_grid2}
P^N(v)=\MT v(x_0^N)&v(x_1^N)&\cdots&v(x_N^N)\EM^T\\
\Phi^N(\MT v_0&v_1&\cdots&v_N\EM^T)=\ds\sum_{k=0}^N v_kq_k^N
\EE
The inner product in $\Real^N$ can be induced from $L^2(\Omega)$, i.e. for $u, v\in \Real^N$, 
$$<u,v>_N=<\ds\sum_{k=0}^N u_kq_k^N, \ds\sum_{k=0}^N v_kq_k^N>$$
If 
$$W=\mbox{span}\{ q^N_0, q^N_1,\cdots, q^N_N\}$$
then the norms $||\cdot ||_N$  and $||\cdot ||_X$ satisfy the consistency assumption (\ref{mapping3b}).$\square$\\

Following \cite{kang-xu1}, we define the observability for ODE systems. 
\begin{definition}
\label{def2}
Given $\rho >0$ and a trajectory $u^N(t)$ of (\ref{eqodemodelnon}) with $u^N(0)\in W^N$. Denote $\hat y^N(t)=\calH\circ\Phi^N(\hat u^N(t))$ and $y^N(t)=\calH\circ\Phi^N(u^N(t))$. Let
$$\epsilon^N = \inf ||\hat y^N-y^N||_Y$$
where $\hat u^N$ satisfies
\EQ
\label{hatuNnon}
d\hat u^N/dt=F^N(t,\hat u^N)\\
\hat u^N(0)\in W^N\\
||\hat u^N(0) -u^N(0)||_N=\rho 
\EE
Then $\rho/\epsilon^N$ is called the unobservability index of $u^N(0)$ in the space $W^N$.
\end{definition}

\section{The consistency of observability}
\label{sec4}
In this section, we prove the consistency of observability. 

\begin{theorem} 
\label{theorem2}
(Consistency) Suppose the initial value problem (\ref{eqpdemodelnon}) and its approximation scheme (\ref{eqodemodelnon})-(\ref{mappingnon}) are well-posed. Suppose $\calH$ satisfies the continuity assumption (\ref{continuityassumption}) in $C^0([0,T],\bigcup_{N=N_0}^\infty \Phi^N(\Real^N))$ for some integer $N_0>0$. Consider a nominal trajectory $u(t)$ of (\ref{eqpdemodelnon}) and its corresponding trajectory $u^N(t)$ of (\ref{eqodemodelnon}) with an initial value $u^N(0)=P^N(u(0))$. Assume that the sphere in $W^N$ centered at $u^N(0)$ with radius $\rho$ is contained in $P^N(D_0)$. Then, 
\EQ
\label{consistnon}
\ds\lim_{N\rightarrow \infty} \epsilon^N =\epsilon
\EE
\end{theorem}

To prove this theorem we need two lemmas. 
\begin{lemma}
\label{lemma1b}
Given any sequence 
$\{ v^N(t)\}_{N=N_0}^\infty$
where $v^N(t)$ is a solution of (\ref{eqodemodelnon}) with $v^N(0)\in P^N(D_0)\bigcap W^N$ and $N_0>0$ is an integer. If $\{ ||v^N(0)||_N\}_{N_0}^\infty$ is bounded and if $\Phi^N(v^N(0))$ converges to $v(0)\in D_0$, where $v(t)$ is a solution of (\ref{eqpdemodelnon}), then $\Phi^N(v^N(t))$ converges to $v(t)$ uniformly for $t\in [0,T]$.
\end{lemma}

{\it Proof }. Let  $\tilde{v}_N(t)$ be the solution of the PDE  (\ref{eqodemodelnon}) such that $\tilde{v}_N(0) \in D_0\bigcap W$ and $P^N(\tilde{v}_N(0))= v^N(0)$. Note that we do not assume $\Phi^N(P^N(\tilde v_N(0)))=\tilde v_N(0)$, although this is the case in many approximation schemes. The set $\{ \tilde v_N(0)\}$ must be bounded in $X$ because of   (\ref{mapping3b}) and the assumption that $\{ ||v_N(0)||_N\}_{N_0}^\infty$ is bounded. Therefore, (\ref{rateconv}) implies 
\EQ
\label{eqtmp_a}
\ds\lim_{N\rightarrow \infty}||\Phi^N(v^N(t))-\tilde v_N(t)||_X=0
\EE
converges uniformly. In particular, 
$$\ds\lim_{N\rightarrow \infty}||\Phi^N(v^N(0))-\tilde v_N(0)||_X=0$$
Because $\Phi^N(v^N(0))$ converges to $v(0)\in D_0$, we know that $\tilde v_N(0)$ converges to $v(0)$. 
Because the solutions of the PDE are continuously dependent on their initial value (well-posedness), 
\EQ
\label{eqtmp_b}
||\tilde v_N(t) - v(t)||_X\rightarrow 0
\EE
uniformly in $t$. Equations (\ref{eqtmp_a}), (\ref{eqtmp_b}) and the triangular inequality
$$||\Phi^N(v^N(t))-v(t)||_X\leq ||\Phi^N(v^N(t))-\tilde v_N(t)||_X+||\tilde v_N(t) - v(t)||_X$$
imply
$$||\Phi^N(v^N(t))-v(t)||_X\rightarrow 0 \mbox{ as } N\rightarrow \infty$$
uniformly in $t$.
$\square$


\begin{lemma}
\label{lemma2}
Given a sequence $\hat u^{N}(t)$, $N\geq N_0$, satisfying (\ref{hatuNnon}). There exists a subsequence, $\hat u^{N_k}(t)$ such that $\{ \Phi^{N_k}(\hat u^{N_k}(t))\}_{k=1}^\infty$ converges uniformly to a solution of (\ref{hatunon}) as $N_k\rightarrow \infty$. 
\end{lemma}

{\it Proof }.  For each $\hat u^N(t)$, there exists a trajectory $\hat u_N(t)$ of the PDE (\ref{eqpdemodelnon}) satisfying $\hat u_N(0)\in D_0\bigcap W$ and $P^N(\hat u_N(0))=\hat u^N(0)$. From the consistency of norms, (\ref{mapping3b}), and the fact $||\hat u^N(0)||_N=\rho$, we know that $\{\hat u_N(t)\}_{N_0}^\infty$ is a bounded set in $W$. The convergence assumption (\ref{rateconv}) implies
\EQ
\label{adeq_c}
\ds\lim_{N\rightarrow \infty}||\hat u_N(0)-\Phi^N(\hat u^N(0))||_X=0
\EE
As a bounded set in a finite dimensional space $W$, $\{\hat u_N(t)\}_{N_0}^\infty$ has a convergent subsequence $\{\hat u_{N_k}(t)\}_{k=1}^\infty$. Let $\hat u(0)$ be its limit and $u(t)$ be the corresponding trajectory of (\ref{eqpdemodelnon}). From (\ref{adeq_c}), we know 
$$\ds\lim_{N\rightarrow \infty}||\Phi^{N_k}(\hat u^{N_k}(0))-\hat u(0)||_X=0$$
From Lemma \ref{lemma1b}, $\{ \Phi^{N_k}(\hat u^{N_k})\}_{k=1}^\infty$ converges to $\hat u(t)$ uniformly in $[0,T]$. Because of (\ref{mapping3b}),
$$\begin{array}{lll}
|| \hat u(0)-u(0)||_X\\
=\ds\lim_{k\rightarrow \infty}||\hat u_{N_k}(0)-u(0)||_X\\
=\ds\lim_{k\rightarrow \infty}(||\hat u^{N_k}(0)-u^N(0)||_N+a_N||\hat u^{N_k}(0)-u^N(0)||_N)\\
=\rho
\end{array}$$
Therefore, $\hat u(t)$ satisfies (\ref{hatunon}) and $\{ \Phi^{N_k}(\hat u^{N_k})\}_{k=1}^\infty$ 
converges uniformly to a solution of (\ref{hatunon}).
$\square$\\

{\it Proof of Theorem \ref{theorem2}}. First, we prove
\EQ
\label{lowlimnon}
\liminf \epsilon^N \geq \epsilon
\EE
Suppose this is not true, then $\liminf \epsilon^N < \epsilon$. There exists $\alpha >0$ and a subsequence $N_k\rightarrow \infty$ so that
$$\epsilon^{N_k} < \epsilon -\alpha$$
for all $N_k$. From the definition of $\epsilon^N$, there exist $\hat u^{N_k}(t)$ satisfying (\ref{hatuNnon}) such that 
$$||\hat y^{N_k}-y^{N_k}||_Y < \epsilon - \alpha$$
where 
$$\hat y^{N_k}(t)=\calH\circ\Phi^{N_k}(\hat u^{N_k}(t)), \;\; y^{N_k}(t)=\calH\circ\Phi^{N_k}(u^{N_k}(t)) $$
From Lemma \ref{lemma2}, we can assume that $\Phi^{N_k}(u^{N_k}(t))$ converges to $\hat u(t)$ uniformly, where $\hat u(t)$ satisfies (\ref{hatunon}). From the continuity of $\calH$ as defined in (\ref{continuityassumption}), 
$$\lim_{k\rightarrow \infty} ||\hat y^{N_k}-y^{N_k}||_Y =||y_{\hat u}-y_u||_Y \leq \epsilon - \alpha$$
However, from the definition of $\epsilon$, we know 
$$ \epsilon \leq ||y_{\hat u}-y_u||_Y$$
A contradiction is found. Therefore, (\ref{lowlimnon}) must hold. 

In the next step, we prove
\EQ
\label{uplimnon}
\limsup \epsilon^N \leq \epsilon
\EE
It is adequate to prove the following statement: for any $\alpha>0$, there exists $N_1>0$ so that 
\EQ
\label{ineqnon}
\epsilon^N < \epsilon +\alpha
\EE
for all $N\geq N_1$. From the definition of $\epsilon$, there exists $\hat u$ satisfying (\ref{hatunon}) so that 
\EQ
\label{eqad1non}
||y_{\hat u}-y_u||_Y < \epsilon +\alpha
\EE
Let $\hat u^N$ be a solution of the ODE
$$du^N/dt=F^N(t,u^N)$$
with an initial value
$$\hat u^N(0)=P^N(\hat u(0))$$ 
Then the following limit converges uniformly
\EQ
\label{eqlimnon}
\ds\lim_{N\rightarrow \infty} ||\Phi^N(\hat u^N(t))-\hat u(t)||_X = 0
\EE
A problem with $\hat u^N(0)$ is that its distance to $u^N(0)$ may not be $\rho$, which is required by (\ref{hatuNnon}). Let $\bar u(t)$ be a solution of (\ref{eqodemodelnon}) with an initial value 
$$\begin{array}{lll}
\bar u^N(0)=\gamma_N(\hat u^N(0)-u^N(0))+u^N(0)\\
\gamma_N=\Fr{\rho}{||\hat u^N(0)-u^N(0)||_N}
\end{array}$$
Obviously 
$$||\bar u^N(0)-u^N(0)||_N=\rho$$
and $\bar u(t)$ satisfies (\ref{hatuNnon}). 
Due to the consistency of the norms and the fact $||\hat u(0)-u(0)||_X=\rho$, we know
\EQ
\label{limgamma}
\lim_{N\rightarrow \infty} \gamma_N=1
\EE
From (\ref{eqlimnon}) and (\ref{limgamma}), we have
$$\begin{array}{lll}
\ds\lim_{N\rightarrow \infty} \Phi^N(\bar u^N(0))\\
=\ds\lim_{N\rightarrow \infty}\left(\gamma_N(\Phi^N(\hat u^N(0))-\Phi^N(u^N(0)))+\Phi^N(u^N(0))\right)\\
=\hat u(0)
\end{array}$$
By Lemma \ref{lemma1b}, the limit 
$$\ds\lim_{N\rightarrow \infty} \Phi^N(\bar u^N(t))=\hat u(t)$$
converges uniformly on the interval $t \in [0,T]$. Let $\bar y^N(t)=\calH(\Phi^N(\bar u^N(t)))$ and $ y^N(t)=\calH(\Phi^N(u^N(t)))$.
Then output continuity assumption (\ref{continuityassumption}) implies
$$\begin{array}{lll}
\ds\lim_{N\rightarrow \infty} ||\bar y^N-y^N||_Y
=||y_{\hat u}-y_u||_Y
<\epsilon + \alpha
\end{array}$$
There exits $N_1>0$ so that
$$||\bar y^N-y^N||_Y < \epsilon+\alpha$$
for all $N\geq N_1$. From the definition of $\epsilon^N$, we know
$$\epsilon^N \leq ||\bar y^N-y^N||_Y< \epsilon + \alpha$$
for all $N>N_1$. Therefore, (\ref{uplimnon}) holds. 

To summarize, the inequalities (\ref{lowlimnon}) and (\ref{uplimnon}) imply 
$$\lim_{N\rightarrow \infty} \epsilon^N =\epsilon$$
$\square$

\section{The empirical Gramian}
\label{sec5}

In this section, we assume that $W^N$ and the space of $y(t)$ are both Hilbert spaces with inner products $<,>_N$ and $<\cdot,\cdot>_Y$, respectively. The ratio, $\epsilon^N/\rho$, is approximately equal to the smallest eigenvalue of a Gramian matrix. More specifically, let 
$$\{e_1^N,e_2^N,\cdots,e_s^N\}$$
be a set of orthonormal basis of $W^N$. Let $u^{+}_i(t)$ and  $u^{-}_i(t)$ be solutions of (\ref{eqodemodelnon}) satisfying
$$u^{\pm}_i(0)=u^N(0)\pm \rho e^N_i$$ 
Define
$$\begin{array}{lll}
\Delta u_i=u^{+}_i(t)-u^{-}_i(t)\\
\Delta y_i(t)=\calH\circ\Phi^N(\Delta u_i(t))
\end{array}$$
The empirical Gramian is defined by
$$\begin{array}{lll}
G=[ G_{ij}]\\
G_{ij}=\Fr{1}{4\rho^2}<\Delta y_i,\Delta y_j>_Y
\end{array}$$
For a small value of $\rho$, the unobservability index is approximated by
\EQ
\label{unobsgramian}
\rho/\epsilon^N\approx\Fr{1}{\sqrt{\sigma_{min}}}
\EE
where $\sigma_{min}$ is the smallest eigenvalue of $G$. For linear ODEs, this approximation is accurate because it can be shown that
$$\begin{array}{rclll}
(\epsilon^N)^2
&=&\ds\min_{\sum a_k^2=\rho^2} \MT a_1&a_2&\cdots&a_s\EM G \MT a_1&a_2&\cdots&a_s\EM^T\\
&=&\sigma_{min}\rho^2
\end{array}$$
Comparing to the linear control theory, $G$ is the same as the observability Gramian if $W^N$ is the entire space and if $y(t)$ lies in a $L^2$-space. 

If $\{e_1^N,e_2^N,\cdots,e_s^N\}$ is not an orthonormal basis, we can modify the approximation as follows. Define
\EQ
\label{eqGU}
S_{ij}=<e_i^N,e_j^N>_N
\EE
Let $\sigma_{min}$ be the smallest eigenvalue of $G$ relative to $S$, i.e.
$$G\xi = \sigma_{min}S\xi$$
for some nonzero $\xi \in \Real^s$. Then
$$\rho/\epsilon^N\approx\Fr{1}{\sqrt{\sigma_{min}}}$$
The approximation is accurate for linear systems.

For the heat equation approximated using (\ref{eqheatode}), the associated mappings can be defined by
$$\begin{array}{lll}
P^N(u)=\MT u^N_1,u^N_2,\cdots,u^N_N\EM^T, &u^N_k=\Fr{2}{L}\ds\int_0^{2\pi} u(x)\sin \left(\Fr{k\pi x}{L}\right)dx\\
\Phi^N(u^N)=\ds\sum_{k=1}^{N} u_k^N\sin \left(\Fr{k\pi x}{L}\right)
\end{array}$$
If we want to find the observability of the first $s$ modes, Definition \ref{def2} is equivalent to the analysis using the traditional observability Gramian for $N=s$. In fact, for all $N\geq s$, $G$ is a constant matrix and 
$$G=\ds\int_0^Te^{(A^s)'t}(C^s)'C^se^{A^st}dt$$
Therefore, $\epsilon^N=\epsilon^s$ for all $N\geq s$ and $\epsilon^N$ is a consistent.

\section{Examples}
\label{sec7}

In this section, some examples are used to illustrate Definition \ref{def1}.
In the literature, a measure of controllability for linear systems is defined based upon the radius of matrices \cite{burns,eising,kenney}. Using duality, we can define the observability radius, $\gamma_o$, as the distance between the system and the set of unobservable systems. This radius is defined primarily for the measure of the robustness of observability, although it agrees with the partial observability defined above in some filtering problems. Using the following example, we illustrate the sameness and differences between $\gamma_o$ and the unobservability index.
\vspace{0.2in}

\noindent {\bf Example 1}. Consider a linear system
\EQ
\label{eq1_ex1}
\MT \dot x_1\\ \dot x_2\EM = \MT 1 & \delta \\ 0 & 1\EM \MT x_1\\ x_2\EM\\
y=x_1
\EE
where $\delta$ is a constant number. In \cite{burns}, it is proved that using Euclidean norm we have
$$\gamma_o=\delta$$
For a small $\delta$, the observability of the system can be qualitatively changed by a small perturbation, i.e. the observability is not robust. This makes sense because, when $\delta=0$, the system is unobservable. Consider the observability of $x(T)$ for some fixed time $T$. The solution of (\ref{eq1_ex1}) is 
\bea
x_1(t)&=&(x_1(T)+\delta x_2(T)(t-T))e^{t-T}\\
x_2(t)&=&x_2(T)e^{t-T}\\
y(t)&=& (x_1(T)+\delta x_2(T)(t-T))e^{t-T}
\eea
It it straightforward to derive the norm of $y$
\bean
\label{eq2_ex1}
||y||^2&=&\int_0^Ty^2(t)dt\\
&=& \MT x_1(T)&x_2(T)\EM \MT \alpha_{11} &\alpha_{12}\delta\\ \alpha_{21}\delta& \alpha_{22}\delta^2\EM \MT x_1(T)\\ x_2(T)\EM
\eean
where 
$$\alpha_{11}=1-e^{-2T},\; \alpha_{12}=\alpha_{21}=(T+\Fr{1}{2})e^{-2T}-\Fr{1}{2}, \alpha_{22}=\Fr{1}{2}-(T^2+T+\Fr{1}{2})e^{-2T}$$
Therefore, the unobservability index of $x(T)$ is
$$\rho/\epsilon=1/\sigma_{min}$$
where $\sigma^2_{min}$ is the smallest eigenvalue of the matrix in (\ref{eq2_ex1}). It can be shown that 
$$\sigma^2_{min}=\Fr{\frac{1}{4}(e^{-4T}+1)-(T^2+\frac{1}{2})e^{-2T}}{1-e^{-2T}}\delta^2+O(\delta^3)$$
If $T$ is large, then 
\EQ
\label{eq3_ex1}
\epsilon/\rho =\sigma_{min}\approx \Fr{1}{2}\delta+O(\delta^3) 
\EE
Because $\gamma_o=\delta$, (\ref{eq3_ex1})  implies that the observability of $x(T)$ defined in this paper coincides with the radius of observability $\gamma_o$. 

Although $\gamma_o$ reflects the observability of the states at $t=T$, $\rho/\epsilon$ and $\gamma_o$ have some fundamental differences. The value of $\gamma_o$ represents the radius of observability which is a measure of observability robustness. It cannot tell the observability of $x(t)$ when $t\neq T$. In fact, for a small value of $\delta$, the initial state $x(0)$ can still be strongly observable if $T$ is large, assuming that $\delta$ is a known constant. Another difference lies in the fact that $\rho/\epsilon$ is defined for the partial observability of a large system, which is not reflected by $\gamma_o$.
$\Box$

In the following example we use Burgers' equation to illustrate the partial observability in a subspace based on the Fourier exansion.
\vspace{0.2in}
 
\begin{figure}[!ht]
	\begin{center}
		\includegraphics[width=3.0in,height=3.0in]{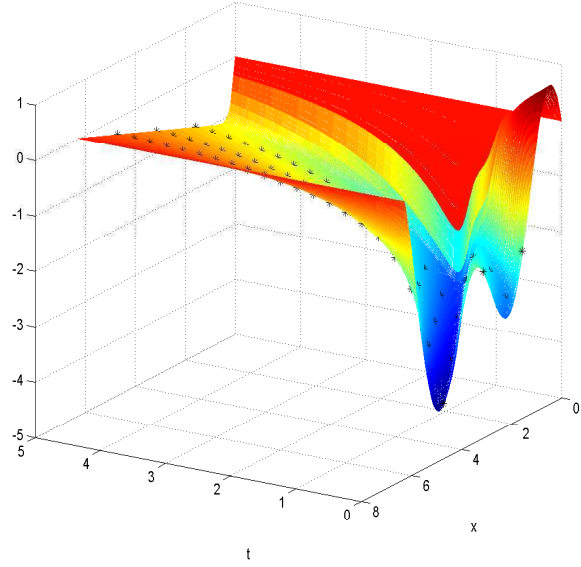} 
		\caption{A solution of Burgers' equation with sensor measurements}
		\label{Figburgers}
		\end{center}
\end{figure}

\noindent {\bf Example 2}. Consider the following Burgers' equation
\EQ
\label{pde}
\Fr{\partial u(x,t)}{\partial t}+u(x,t)\Fr{\partial u(x,t)}{\partial x}=\kappa \Fr{\partial^2 u(x,t)}{\partial x^2} \\
u(x,0)=u_0(x), & x\in [0, L]\\
u(0,t)=u(L,t)=0, & t\in [0,T]
\EE
where $L=2\pi$, $T=5$, and $\kappa=0.14$. The output space consists of functions representing data from three sensors that measure the value of $u(x,t)$ at fixed locations 
\EQ
\label{pdeoutput}
\MT y_1(t_k) \\ y_2(t_k)\\ y_3(t_k)\EM=\MT u(\frac{L}{4},t_k)\\ u(\frac{2L}{4},t_k)\\ u(\frac{3L}{4},t_k)\EM, & t_k=k\Delta t, \;\; k=0,1,2,\cdots,N_t
\EE
where $\Delta t=T/N_t$, $N_t=20$. Figure \ref{Figburgers} shows a solution with discrete sensor measurements marked by the stars. In the output space
$$||y||_Y=\left(\dsum_{k=0}^{N_t}(y_1^2(k)+y_2^2(k)+y_3^2(k))\right)^{1/2}$$

The approximation scheme is based on equally spaced grid-points
$$x_0=0 < x_1 <\cdots < x_N=L,$$ 
where 
$$\Delta x=x_{i+1}-x_i=L/N.$$  
System (\ref{pde}) is discretized using a central difference method
\EQ
\label{odemodel}
\dot u_1^N=-u_1^N\Fr{u_2^N-u_0^N}{2\Delta x}+\kappa \Fr{u_{2}^N-2u_1^N}{\Delta x^2}\\
\dot u_2^N=-u_2^N\Fr{u_3^N-u_1^N}{2\Delta x}+\kappa \Fr{u_{3}^N+u_{1}^N-2u_2^N}{\Delta x^2}\\
\;\;\;\;\;\vdots\\
\dot u_{N-1}^N=-u_{N-1}^N\Fr{u^N_{N}-u_{N-2}^N}{2\Delta x}+\kappa \Fr{u_N^N+u_{N-2}^N-2u_{N-1}^N}{\Delta x^2}
\EE
where $u^N_0=u^N_{N}=0$. For any $v(x)\in C([0,L])$, we define
$$P^N(v)=\MT v(x_1) & v(x_2) & \cdots v(x_{N-1})\EM \in \Real^{N-1}$$
For any $v^N\in \Real^{N-1}$, define 
$$\Phi^N(v^N)=v(x)\in C[0, L]$$ 
be the unique function of cubic spline determined by $v^N$ and $(x_0,x_1,\cdots,x_N)$ satisfying $v(0)=v(L)=0$. We adopt $L^2$-norm in $C[0, L]$. For any vector $v^N \in \Real^{N-1}$, its norm is defined as 
$$ ||v^N||^2_N= \Fr{2\pi}{N} \sum_{i=1}^{N-1} v^2_i$$

The space for estimation is defined to be 
$$W=\left\{\alpha_0/2+ \dsum_{k=1}^{K_F} \left(\alpha_k\cos(\Fr{2k\pi}{L}x)+\beta_k \sin(\Fr{2k\pi}{L}x)\right)\left |\begin{array}{lll}\alpha_k, \beta_k \in \Real\\ \alpha_0/2+ \dsum_{k=1}^{K_F} \alpha_k=0\end{array}\right. \right\} $$
In this section, $K_F=2$.  This means that we want to find the observability for the first five modes in the Fourier expansion of $u(0)$. Or equivalently, we would like to find the observability of 
$$\MT\alpha_0 & \alpha_1 &\beta_1 &\alpha_2 &\beta_2\EM$$
Define 
$$X^N=\MT x_1 & x_2&\cdots &x_{N-1}\EM^T$$
then
$$W^N=\left\{\alpha_0/2+ \dsum_{k=1}^{K_F} \left(\alpha_k\cos(\Fr{2k\pi}{L}X^N)+\beta_k \sin(\Fr{2k\pi}{L}X^N)\right)\left |\;  \begin{array}{lll}\alpha_k, \beta_k \in \Real\\ \alpha_0/2+ \dsum_{k=1}^{K_F} \alpha_k=0\end{array}\right. \right\} $$

\begin{figure}[!ht]
	\begin{center}
		\includegraphics[width=4.0in,height=2.0in]{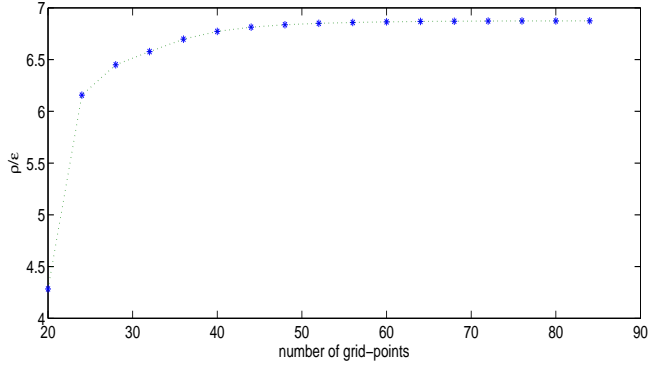} 
		\caption{The observability of Burgers' equation}
		\label{figburger}
		\end{center}
\end{figure}

In this example, the nominal trajectory has the following initial value 
$$u_0(x)=-2+\cos(x)+\sin(x)+\cos(2x)+\sin(2x)$$
Its solution is shown in Figure \ref{Figburgers}. To approximate its observability, we apply the empirical Gramian method to (\ref{odemodel}) in the space $W^N$. The ratio $\rho/\epsilon^N$ is approximated for $N=4k$, $5\leq k \leq 21$. The value of unobservability index approaches (Figure \ref{figburger}) $$\rho/\epsilon=6.87$$ 
$\Box$ 

The partial observability is defined for PDEs approximated using general approximation schemes. Different from the finite difference method adopted in Example 2, in Example 3 the partial observability is computed for a shallow water model derived by using a finite element method. 
\vspace{0.2in}

\begin{figure}[!ht]
	\begin{center}
		\includegraphics[width=4.0in,height=2.0in]{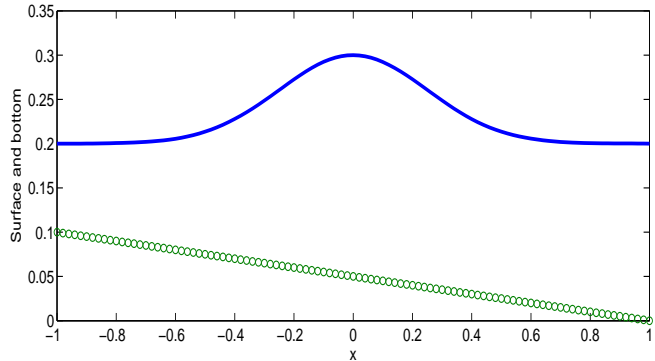} 
		\caption{The initial surface height (solid line) and the bottom topography height (circle)}
		\label{figSW_initial}
		\end{center}
\end{figure}

\noindent {\bf Example 3}. The shallow water equations are widely used in scientific computation, numerical weather prediction, and oceanography as a testbed or illustrative example. For the one dimensional case, the independent variables include time, $t$, and the space coordinate, $x$. The dependent variables are the fluid depth, h(x,t), and the horizontal velocity $u(x,t)$. The differential equations are
\bean
\label{eq1_ex3}
&&\Fr{\partial h}{\partial t}+\Fr{\partial (uh)}{\partial x}=0\\
&&\Fr{\partial (uh)}{\partial t}+\Fr{\partial (u^2h+\frac{1}{2}gh^2)}{\partial x}=0
\eean
Assume that the sensors are located at $x=0.2$, $0.5$, and $0.8$. More specifically
\EQ
\label{SWoutput}
\MT y_1(t_k) \\ y_2(t_k)\\ y_3(t_k)\EM=\MT h(0.2,t_k)\\ h(0.5,t_k)\\ h(0.8,t_k)\EM, & t_k=k\Delta t, \;\; k=0,1,2,\cdots,N_t
\EE
where $\Delta t=1.25\times 10^{-3}$ and $N_t=800$. The initial surface height is defined by a Gaussian curve (Figure \ref{figSW_initial})
$$h_0(x)=0.1\exp \left(-8(x-\frac{1}{2})\right)+0.2$$
where $x\in [-1, 1]$ and $t\in [0,1]$. The height of the bottom topography is a line  (Figure \ref{figSW_initial})
$$h_b(x)=0.1(1-x)$$ 

\begin{figure}[!ht]
	\begin{center}
		\includegraphics[width=4.0in,height=2.0in]{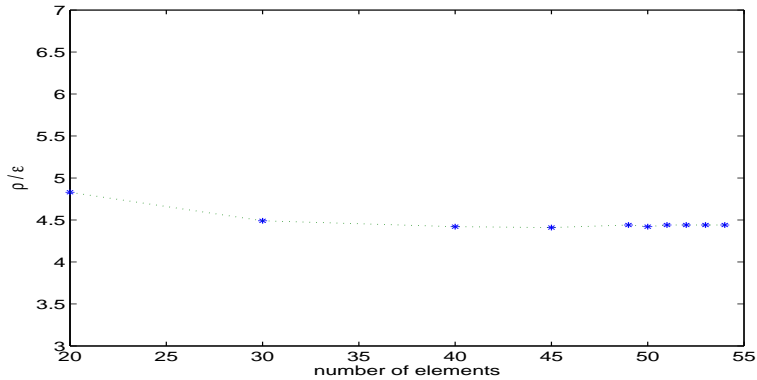} 
		\caption{The observability consistency of the Shallow Water Model}
		\label{figSW}
		\end{center}
\end{figure}

Based on the empirical Gramian method, the partial observability can be approximated using a computer code that generates the states of (\ref{eq1_ex3}). The original differential equation is not directly used in the computation. In this example, we adopt a code that discretizes (\ref{eq1_ex3}) using a spectral element method. More specifically, the space $ [-1, 1]$ is divided into finite many subintervals or elements, $I_1$, $I_2$, $\cdots$, $I_N$. In each element $I_i$, the solution is approximated using Lagrange polynomials at the Legendre-Gauss-Lobatto (LGL) nodes, $x_{ij}$, $j=0,1,\cdots,N_p$, where $N_p$ is the order of the polynomials on the interval $I_i$. The approximate solution has the form
$$h_{ij}(t_k)\approx h(x_{ij},t_k), \;\; u_{ij}(t_k)\approx u(x_{ij},t_k)$$
For the mapping $\Phi^N$, we simply use the linear interpolation. Note that this is different from the spectral element method in the discretization of the PDE. In fact, in the computation of observability the approximation scheme can be different from the method used in the PDE discretization. Suppose the initial height, $h(x,0)$, has a Fourier expansion. The space for estimation is 
$$W=\left\{\alpha_0/2+ \dsum_{k=1}^{K_F} \left(\alpha_k\cos(\Fr{2k\pi}{L}x)+\beta_k \sin(\Fr{2k\pi}{L}x)\right)\left |\alpha_k, \beta_k \in \Real\right. \right\} $$
where $K_F=6$. We would like to compute the observability of $\alpha_k$, $k=0,1,\cdots, 6$ and $\beta_k$, $k=1,2,\cdots,6$. Under the $L_2$ norm for both the initial $h(x,0)$ and the output function, the unobservability index is shown in Figure \ref{figSW}. After the number of elements is increased to $N=54$, the index is stabilized at $\rho/\epsilon =4.44$.

\section{Conclusions}
A definition of partial observability using dynamic optimization is introduced for PDEs. The advantage of this definition is to resolve several issues and concerns about observability in a unified framework. More specifically, using the concept one can achieve a quantitative measure of partial observability for PDEs. Furthermore, the observability can be numerically approximated. A practical feature of this definition for infinite dimensional systems is that the observability can be numerically computed using well-posed approximation schemes. It is mathematically proved that the approximated partial observability is consistent with that of the original PDE. A first order approximation is derived using empirical Gramian matrices. The concept is illustrated using several examples.  \\
\\
{\bf Acknowledgement}:  The authors would like to express their gratitude to Professor Arthur J. Krener for his comments and discussions on the observability Gramian and the unobservability index. We would like to extend our thanks and appreciation to Professor Tong Huang for his suggestions on some basic assumptions made for the PDE and its solutions.

\end{document}